\documentclass[a4paper,10pt]{article}
\usepackage{amssymb}
\usepackage[T2A]{fontenc}
\usepackage[cp1251]{inputenc}
\usepackage{amsmath}

\usepackage[dvips]{graphicx}
\usepackage[english,russian]{babel}

\newtheorem{theorem}{Теорема}[section]

\newtheorem{definitionhead}[theorem]{Определение}

\author{А.~Я.~Белов}

\title{Олимпиады: дверь в математику или спорт?}

\begin{document}

\maketitle

\section{Введение}

По всему миру проводятся математические конкурсы и олимпиады.
Появились специалисты по их проведению, возникла олимпиадная
математика со своей методикой работы и своей литературой.
Олимпиадный мир стал жить собственной жизнью, но его кажущаяся
самодостаточность породила ряд проблем.

О пользе и вреде олимпиад, о том, как их проводить, постоянно
ведутся кулуарные дискуссии. Во время заседаний методических
комиссий, когда следует принимать решения, эти дискуссии
превращаются в довольно трудные разговоры. Поэтому необходимо их
провести открыто на страницах журнала. Сожалея о неизбежной субъективности,
с благодарностью выслушаю замечания.

Ситуация вокруг олимпиад  парадоксальна. С одной стороны, на
олимпиады тратятся значительные ресурсы. Талантливых школьников и
их учителей сводят вместе, прежде всего, математические олимпиады.
У истоков олимпиадного движения стояли великие ученые. Многие
школьники, особенно на периферии, получают математическое
образование, нацеленное в первую очередь на подготовку к
олимпиадам, со всеми его плюсами и минусами. (С этим надо
считаться научным руководителям и организаторам учебного
процесса.)

С другой стороны, распространены суждения о вреде олимпиад,
зачастую весьма странные. Например, бытует мнение о том, что
успехи на олимпиадах мало связаны с научной карьерой. Доходит до
курьезов~-- до утверждений о том, что ``среди крупных ученых нет
победителей олимпиад'', хотя среди бывших победителей олимпиад
известных математиков во много раз больше, чем среди
неолимпиадников. Автор сталкивался с попытками одного
преподавателя забрать своих учеников из его кружка. Один из
аргументов, высказанных им школьникам, звучал так: ``Концевич в
олимпиадах не участвовал''.

Такого рода суждения очень легко критиковать, а иногда~-- и
высмеивать. Владимир Соловьев высказал мысль о том, что любая
ложная социальная теория базируется на искажении некоторой правды,
и для победы над ``теорией'' необходимо эту внутреннюю правду
осознать и выявить. Мне представляется, критики олимпиад чувствуют
серьезные проблемы, зачастую не умея их
сформулировать.\footnote{Есть математики-спортсмены, нацеленные на
решение конкретных проблем, и математики-``домостроители'',
занимающиеся построением теорий, они чаще всего и становятся
критиками олимпиад.}

Хотя олимпиадное движение играет значительную роль в
математическом образовании школьников, особенно на периферии,
олимпиадные деятели зачастую заявляют, что цель олимпиад
ограничивается выявлением талантливых учащихся и формированием
первоначального интереса к предмету. При этом декларируют, что
``олимпиада~-- это не математика'' (а то, что преподают в школе
или колледже, это математика?), и любят подчеркивать, что далеко
не все становятся математиками (хотя содержание олимпиадных задач
должно быть полноценным, вне зависимости от будущей профессии).
Вовлекать в математику надо с помощью математики же, а выявлять
таланты -- с помощью полноценного материала. \footnote{Некоторые
чрезвычайно важные вещи могут делаться только в качестве гарнира,
хотя этот гарнир зачастую бывает полезнее основного блюда.}

В адрес некоторых задач со стороны некоторых активных деятелей
олимпиадного движения высказывалась и такая критика, пусть и в
полемическом запале: {\it Плохо, что эта задача из науки}.
Постоянно ведутся разговоры о спортивных достижениях. Победа
школьника на олимпиаде иногда ценится выше его публикации в
академическом журнале.

Против тезиса о важности научного содержания олимпиад
высказывается и такой странный аргумент: ``Олимпиады~-- это только
небольшой жизненный эпизод''. Даже если считать прямое и косвенное
действие олимпиадного мира на подростка незначительным (что
неправда )\footnote{На самом деле подготовка и участие в
олимпиадах требует значительных ресурсов. То, с чем знакомиться
ученик в этом процессе отпечатывается на всей его будущей карьере.
Влияет все это и на учителей.} , этот аргумент выглядит странно.
Следуя такой логике, можно оправдать плохое качество любого
отдельно взятого урока, ибо один урок мало что решает. Эту
аргументацию не следует воспринимать как аргументы в научном
споре, но выявлять скрытые мотивы -- не спорить, но объяснять. Я
ни разу не слышал ничего подобного со стороны действующих ученых.
Уход от ответственности за научное будущее ученика узким
олимпиадным деятелям необходим, чтобы оправдать автономное
существование олимпиадного мира.

Работа над этой статьей началась еще в 1996 году. В дальнейшем
стало понятно, что проблемы носят не только внутрироссийский, но и
международный характер. Так, участники команды одной страны, не
набравшие достаточного количества баллов, не были допущены до
фотографирования. На международных  олимпиадах по математике
неоднократно наблюдались случаи коррупции (сообщение задач членам
своей команды или ``пробивание'' в вариант задач, которые команда
знает, пристрастная проверка работ, сбор компромата на участников
других команд и т.д.), от чего сильно страдают и содержание
вариантов, и результаты участников. Разговоры педагогов о
школьниках напоминают разговоры о скаковых лошадях. Другой
побудительной причиной написания этой статьи послужили внутренние
проблемы, связанные с самим олимпиадно-педагогическим сообществом.

Данная статья призвана инициировать постановку и обсуждение
назревших проблем, поэтому количество затронутого материала
выходит за формат обычной статьи. Многие аспекты, как нам кажется,
ранее не только не рассматривались должным образом, но даже и не
затрагивались. Каждому из них следовало бы посвятить отдельную
работу.

\paragraph{Благодарности.}  Автор признателен
А.~И.~Буфетову, Э.~Б.~Винбергу,
 М.~Н.~Вялому, А.~Домошницкому,
А.~К.~Ковальджи,
 В.~Н.~Латышеву,
Н.~Х.~Розову,
 В.~М.~Тихомирову,
Б.~Р.~Френкину, Г.~Гусеву, А.~К.~Кулыгину за полезные обсуждения и
поддержку. Автор признателен своим коллегам по проведению
олимпиад, благодаря существованию  которых  эта статья была
написана.

\paragraph{Зачем нужны олимпиады?}
Предваряя обсуждение, автор считает необходимым обозначить личное
отношение к олимпиадам и олимпиадной математике. К сожалению,
обоснование этой точки зрения выходит за рамки настоящей статьи.
Автору близки позиции, изложенные в работах
\cite{Skopenkov1,konstantinov2}.

Решение олимпиадных задач (разного уровня сложности) служит
основой для почти всех математических кружков. Подготовка к
олимпиадам оказывает значительное влияние на первоначальные
занятия школьников математикой. Именно в решении трудной задачи
может состоять достижение подростка. Более того, он зачастую
оказывается в равном положении со взрослым. На трудных задачах
вырабатывается интеллектуальная техника и соответствующие волевые
качества. Но главное~-- сам факт достижения серьезной, но
посильной цели в подростковом возрасте.

Часто утверждается, что олимпиадные умения не связаны с большой
наукой, а зависимость между олимпиадными успехами и научной
карьерой весьма слабая. Автор с этим категорически не согласен,
поскольку препятствий для развития таланта множество. Прежде
всего, человек даже очень талантливый встречается с теми или иными
житейскими обстоятельствами. Они его могут надломить или даже
сломать. Он может уйти в зарабатывание денег, столкнуться с
семейными проблемами. Поэтому как бы мы ни выявляли таланты в юном
возрасте, какая-то часть (и, увы, очень большая) из них в зрелом
возрасте погаснет. (Небольшая популяция в Древней Греции поставила
много великих ученых~-- больше, чем та же Греция произвела за
последние две тысячи лет). Не исключено, что олимпиадные успехи
больше говорят об изначальном таланте, чем будущее научное
творчество.

Деятельность вокруг олимпиад стала заметным явлением в области
современного математического образования. Их роль далеко не
ограничивается обнаружением талантливых учащихся. Благодаря
олимпиадной математике удается увидеть роль стандартных идей и
рассуждений. Появились подборки олимпиадных задач по темам
``принцип\linebreak  Дирихле'', ``правило крайнего'',
``инварианты'' и др., объединенные единством метода.

Математики-непедагоги, как правило, не уделяют должное внимание
``тривиальным'' вещам, которые между тем играют исключительную
роль как в мышлении математика, так и в подготовке к олимпиадам.

Благодаря олимпиадам возникли знаменитые книги Д.~Пойа
\cite{Polya1, Polya2, Polya3}. Да и облик так называемой
``венгерской математики'' (вспомним Пола Эрдёша) сформировался во
многом под влиянием олимпиад.

Возможно, что выделение стандартных рассуждений может привести к
революции и в самой математике, а в дальнейшем~-- и физике. Сама
работа над изложением, казалось бы, известных результатов, часто
приводила к открытиям. Так было со схемами Дынкина и с уравнением
Гейзенберга. Олимпиадная математика с ее систематизацией идей и
методов может послужить детонатором. Процесс детонации может
начаться в комбинаторике. Мне представляется, что эта наука должна
быть организована не так, как классическая область математики, а
подобно некоторым тематическим подборкам олимпиадных задач. В
основе ее организации должно лежать единство метода. Правильный
учебник по комбинаторике должен быть чем-то вроде учебника
шахматной игры или олимпиадного са\-моучителя.

\paragraph{Из истории олимпиад.}
Олимпиадное движение возникло свыше ста лет тому назад. Первые
олимпиады состоялись в 1884 году в Австро-Венгрии. Они возникли из
конкурсных экзаменов. Затем олимпиады появились в Венгрии. В
дальнейшем олимпиадное движение ширилось. Олимпиады вышли за рамки
конкурсных задач. В 30-е годы по инициативе Б.~Н.~Делоне возникли
Ленинградские, а затем~-- Московские городские олимпиады (Первая в
СССР математическая олимпиада для школьников была проведена в
Тбилиси в декабре 1933 г.). У истоков олимпиад в СССР стояли
ведущие ученые~-- А.~Н.~Колмогоров, И.~М.~Гельфанд,
П.~С.~Александров, С.~Л.~Соболев, Л.~Г.~Шнирельман и другие. В
дальнейшем возникла целая система национальных и международных
олимпиад. Олимпиадный мир стал жить собственной жизнью. До
недавнего времени его лидерами были хорошие математики
 (в~том числе и создатель
питерской олимпиадной школы, которую сейчас некоторые критикуют за
излишне спортивную направленность). Не все знают, что создатель
Турнира городов Н.~Н.~Константинов  имеет красивые математические
результаты.

Подробнее об истории математических олимпиад см. предисловия к
книгам \cite{Leman,Moscow}, а также размышления В.~М.~Тихомирова в
книге \cite{MMO932005}.

\paragraph{Кто сейчас делает олимпиады?}
Хотя у истоков олимпиад стояли великие ученые, в последующем
уровень олимпиадных деятелей постепенно снижался. (Затем~-- раскол
научного сообщества, вызванный проблемами 70-х и начала 80-х
годов, события 90-х годов  усилили эту тенденцию.) Ослабла связь
олимпиад с научным сообществом. Появились так называемые
``олимпиадные функционеры'', ~--е специалисты по организации и
проведению олимпиад. В жюри многих турниров высокого уровня почти
не осталось профессиональных математиков даже среднего уровня.
В~последние годы в некоторых странах появились олимпиадные
деятели, представляющие собой ``нематематиков'' и в то же время
пытающиеся доминировать в олимпиадном мире, иногда откровенно
противопоставляя себя научному сообществу.

Возникшая проблема является относительно новой. Как мне кажется,
она в значительно меньшей мере наблюдается в олимпиадах по другим
предметам, в силу меньшей развитости соответствующих субкультур.
Чтобы ``олимпиадная математика'' смогла сложиться, был необходим
первоначальный приток идей из большой науки. Сейчас, однако, такая
необходимость многими олимпиадными деятелями не только не
ощущается, но ими оказывается сопротивление переносу идей из этого
источника.

\subsection{О вкусах спорят!}
На это мне указал замечательный математик и человек, ныне
покойный, Ромен Васильевич Плыкин (он был организатором и
председателем жюри Всероссийских конференций школьников).
Очевидно, что есть ``прекрасное'' и ``безобразное'' в жизни,  есть
понятие плохого и хорошего вкуса в живописи, в одежде и в еде.
Можно иметь предпочтения в живописи или предпочитать китайскую
кухню французской. Но если вкус человека, разбирающегося во
французской кухне, вине или китайском чае, заслуживает уважения,
то вкус любителя фастфуда~-- нет (потребитель в обществе
потребления и потреблять-то не умеет!).

Любая сильная идея является в обличии красоты (beauty is power
itself). Значение математика определяется не только его
``пробивной силой'' (~--е возможностью ``пробить'' трудную
задачу), но и вкусом, которые, впрочем, тесно связаны. За свой
вкус математик отвечает своей судьбой. Плохая эстетика задач
наносит ущерб учащимся.

На мой взгляд, имеется группировка вкусовых предпочтений членов
жюри олимпиад в зависимости от того, являются ли они действующими
учеными или только функционерами. Если этот факт получит
дополнительное подтверждение, то о плохих или хороших вкусах можно
будет говорить более объективно.

В свое время за счет математического профессионализма и, как
следствие, лучших эстетических критериев подбора задач, в одном
областном центре удалось добиться результатов лучше, чем в
мегаполисах.
(В дальнейшем финансирование талантливых детей в этом
регионе пошло по сомнительному направлению,  в том числе
математические лагеря перестали проводиться.)

\section{Два подхода к олимпиадам}

В олимпиадном мире сложились две ценностные ориентации. Они
проявляются во всем: в подборе задач, выработке критериев оценок,
и~-- что немаловажно~--  отражаются на роли и авторитете тех или
иных личностей и, как следствие, на кадровых вопросах. Проявляются
эти подходы и в организации математических лагерей. В значительной
степени люди привержены одному из этих  подходов, причем не всегда
осознанно. Поэтому автор пытается дать описание, условно говоря,
``научного'' и ``спортивного'' стиля олимпиад.

\subsection{Перерождение в большой спорт}

В последнее время в олимпиадном мире усилился и доминирует
``спортивный'' подход. Я получил публичный упрек от одного
олимпиадного деятеля, что для меня ``{\it олимпиада лишь средство
обучения математике\textup, а для него~-- СПОРТИВНОЕ
СОРЕВНОВАНИЕ}''.

Распространение подобной точки зрения связано с тремя группами
причин. Во-первых, с общим спортивным духом нашего времени.
Во-вторых, с недостатком действующих математиков в олимпиадном
мире и, как следствие, с ухудшением качества кадров. И в-третьих,
это следствие логики развития олимпиадного движения без обратных
связей, которой надо противостоять. (Иногда в предметах с менее
долгой олимпиадной традицией научный уровень жюри бывает выше,
поскольку первоначальный импульс в них задают крупные ученые, по
той же причине бывает выше научный уровень жюри в странах с более
молодой олимпиадной традицией). Ситуация усугубляется влиянием
бизнеса, который в ряде случаев извращает творческие конкурсы (не
только по математике).

Спортивный подход выражается фразой: ``олимпиада~-- это спорт по
решению головоломок''. Этот лозунг влечет за собой многое:
усиливается тренерство, ужесточаются формальные требования и,
соответственно, критерии оценок. Так, если спортсмен переступит
черту на 10 см, а \mbox{прыгнет} на 10 метров, ему не зачтут
прыжок 9,9 м, его прыжок не зачтут вовсе.

Показательны слова одного из олимпиадных деятелей, сказанные во
время обсуждения описки учащегося: ``А если тебе зарплату не так
подсчитали?'' Налицо непонимание цели и смысла олимпиад.

К придумыванию новых задач относятся как к составлению шахматных
этюдов и головоломок, только вместо фигур комбинируют объекты из
школьной программы и стандартные олимпиадные темы. А в шахматах
вопроса ``откуда такое расположение фигур'' просто не возникает.
Комбинировать  пытаются всё со всем. Особенно ценится внешняя
обертка.  Отношение к задачам как к головоломкам ведет к
возникновению химер~-- когда комбинируется несовместимое по своей
внутренней природе. Вот типичные примеры такого рода задач:

\smallskip

{\bf 1.} {\it Можно ли расставить числа от $1$ до $100$ в ряд
так\textup, чтобы сумма любых трех, идущих подряд, была простым
числом?}

{\bf 2.} {\it Стороны треугольника~-- простые числа. Может ли его
площадь быть целым числом?}

\smallskip

В этих двух задачах простота числа притянута искусственно.
Используется только свойство нечетности.

Или еще~--

\smallskip

{\bf 3.} {\it Найдите все целые числа, равные сумме факториалов
своих цифр.}

\smallskip

Это мертвая математика.

\paragraph{Примеры отторжения содержательных задач (разных авторов).}
Речь пойдет о вкусе жюри, а не об аргументации, связанной с
известностью или уровнем сложности задач.

Следующие две задачи были сочтены методической комиссией
``неестественными''.

\smallskip

\textbf{1.} {\it Поезд ехал один час от пункта $A$ в пункт
$B$\textup, проехав $60$~км. Доказать\textup, что в какой-то
момент его ускорение было не менее $240\;
\text{км}/\text{ч}^2$\/}. (В.~М.~Тихомиров)

\textbf{2.} {\it Плоскость покрыта единичными кругами. Докажите,
что некоторая точка покрыта не менее трех раз.} (А.~Я.~Белов)

\smallskip

Вторая задача отражает в простейшей форме фундаментальное понятие
{\it топологической размерности\/}. Пространство имеет {\it
размерность $n$}, когда имеются сколь угодно мелкие покрытия без
перекрытий по $n+2$, но нельзя избавиться от перекрытий по $n+1$.
У того, кому  она кажется неестественной, скорее всего, плохой
вкус и он плохо знает математику.

На одном фестивале была предложена задача по нахождению угла между
некоторыми диагоналями правильного додекаэдра (Автор --
С.~Анисов). Идея решения состояла в рассмотрении вписанного куба.
Эта задача была отвергнута как ``неолимпиадная''. Однако в
задачных конкурсах до недавнего времени использовались разного
рода стереометрические ``монстры''.  Способный от природы школьник
всё же может увидеть куб, вписанный в додекаэдр, а натасканный на
``стандартные'' олимпиадные темы ``спортсмен'', скорее всего, не
увидит.

С другой стороны, многие задачи, например, на построение
инвариантов, могут быть решены только учащимися, хорошо владеющими
этой техникой (которая, к сожалению, даже намеками не входит в
школьный курс). Здесь возникает проблема ``джентльменского
набора'' идей и методов, без владения которыми ``самородок'' не
достигнет больших успехов на олимпиаде. Получается, что следование
вкуса жюри дает преимущество школьников натасканных на олимпиады
по сравнению с изначально талантливыми.

Еще пример ``неолимпиадной задачи''.

{\it Ломаная делит круг на две равные части. Доказать, что она
проходит через его центр.} (А.~К.~Ковальджи).

Стиль решения  этой задачи непривычен для олимпиадных деятелей.
Тут дело вовсе не в трюке.  Надо осознать, что значит {\it две
равные части}. Это значит, что {\it есть движение, переводящее
одну часть в другую\/}. А все типы движений плоскости описаны в
теореме Шаля. Далее следует небольшой перебор. (Подробнее~-- см.
``Математическое просвещение'', сер.~3, вып.~6, 2002.
С.~139--140.)

Или еще пример стереометрической задачи, отторжение которой
говорит о дурном вкусе.

{\it Можно ли разбить пространство на усеченные октаэдры?}

Разговор о конкретных вариантах олимпиад, плохих и хороших задачах
давно назрел. К сожалению, здесь мы имеем возможность только
поставить вопрос о его необходимости.

\paragraph{Причины появления задач-химер и отторжения
содержательных задач.}

Почему получили распространение задачи-химеры? И почему спортивный
подход приводит к негативному, при прочих равных, отношению к
задачам, за которыми стоит внутреннее содержание~-- т.е. поле идей
и сюжетов (т.е. теми, на которых и надо учить математика)?

Дело в том, что поучительная задача чему-то учит. Но тогда и
обратно~-- решение такой задачи непредсказуемым образом зависит от
особенностей решателя и его культуры. Если процесс решения задачи
оказывает воздействие на культуру решателя, то его результаты, в
свою очередь, должны от этой культуры зависеть. Эта зависимость
тем менее предсказуема, чем более глубокой оказывается задача.
Следовательно, такая задача неудобна в плане оценки ее сложности.
Например, нашедшему 100-ю цифру после запятой числа
$(\sqrt{5}+2)^{1000}$ гораздо легче показать несуществование
рациональных $\alpha_i,\beta_i$ таких, что $\sum_i
(\alpha_i+\sqrt{3}\beta_i)^2=2+3\sqrt{3}$. И наоборот, решение
последней задачи облегчает первую. Кроме того, доказавший
разложимость многочлена с вещественными коэффициентами в
произведение многочленов не более чем второй степени имеет
преимущество в решении этих задач. Сейчас, правда, идея
сопряженности обсуждается на занятиях кружка и для прошедшего эти
занятия упомянутые выше задачи суть упражнения. Тем не менее так
было не всегда. Есть много идей и связей которые мы пока не
понимаем. \footnote{Мы предполагаем в дальнейшем добавить
проиллюстрации.}

Спортивный принцип предполагает стандартизацию. Одна из его
основ~-- использование относительно стандартных приемов решения
задач. При решении искусственных задач участники более равны, а
самые ``равные'' должны получать премии.

Большой спорт тяготеет к ограничению поля деятельности и четкой
формализации правил. Поэтому не случайна узость тематики задач,
отсюда опасность вырождения олимпиад.
Кроме того, стиль решения содержательной задачи (за которой стоит
целое поле идей и сюжетов) непривычен для нематематика, а
следовательно, не соответствует его вкусам.


\paragraph{Манипуляторство при проведении занятий.}
Разница между содержательным и спортивным стилем олимпиад примерно
такая же, как  между задачами придумать способ сборки кубика
Рубика и соревнованиями на скорость его сборки. Спортивный стиль,
``головоломочность'' оказывают влияние и на ведение занятий.
Внимание смещается с внутренней сущности на формальные манипуляции
материалом.

Применительно к преподаванию это приводит не только к упору на
натаскивание к олимпиадам. Парадоксальным образом зачастую
наблюдается любовь тренеров к ученым словам и манипулирование ими.
(Примеры курьезов такого рода~-- мини-курс на тему: ``три
определения комплексного числа с доказательством их
равносильности''\footnote{Вполне осмысленно анализировать разные
определения выпуклости фигур, поскольку это дается сравнительно
легко и помогает решать задачи. А в ``игре'' с комплексными
числами, с одной стороны, имеется стремление подражать ``большой
науке'', а с другой~-- мало содержания.}, абстрактно изучаются
``$n$-арные операции'', обсуждаются ``результаты'' типа такого:
группа $S_3$ вкладывается в группу автоморфизмов свободной группы
с тремя образующими и т.д. и т.п.)

На наш взгляд, важна связь учителя с живым источником, которая
сама по себе служит опорой и дисциплинирующим началом, что
позволяет быть менее формальным. Жесткость и формальность в
преподавании связана и с узостью кругозора (чем уже, тем жестче).
Отсутствие или слабость живой связи с наукой приводит к
возрастанию роли внешней обертки, а сами математические понятия
становятся чем-то вроде заклинаний.

\paragraph{Причины распространения формально-спортивного подхода.}

Помимо уже обсуждавшейся логики организации соревнований в большом
спорте есть и другая столь же важная причина распространения
сверхспортивного стиля. Дело в том, что взгляд на математику как
на науку о решении занимательных задач и головоломок самый
доступный. Более того, он необходим при первоначальном знакомстве
с математикой, а следовательно, и в преподавании. Естественно, что
этот самый доступный, безусловно, ценный и живой взгляд на
математику, при всей его узости, получил распространение. Но при
доведении его до крайности возникает сверхспортивный стиль.

Глубина понимания без узости объекта изучения сразу не достигается
(с этим связан подростковый экстремизм типа ``ничего мне не нужно,
кроме геометрии''). Лучше вначале достичь глубины, чем широты. В
доступности спортивного стиля есть и позитивная сторона.
Олимпиадный тренер даже спортивного толка может много сделать для
развития образования в своем регионе.

Талантливый школьный учитель или местный деятель образования  (а
то обстоятельство, что он смог возвыситься над рутиной, говорит о
многом), совершив усилие, иногда даже сверхусилие, входит в
олимпиадный мир. Но чтобы ему понять, что мотивировки задач
принципиально важны, требуется еще одно усилие, которое редко
когда совершается. Помимо всего прочего, человек горд собой~-- у
него уже есть результаты, а они зачастую ослепляют (вплоть до
снобизма). \footnote{Другой механизм пополнения кадров -- через
студентов. При этом  важно, чтобы человек продолжал заниматься
наукой либо хотя бы ее популяризацией и не скатился в чистую
``олимпиадчину'' (а такая опасность есть).} )

\paragraph{Социальные опасности.}

Олимпиадный мир~-- это только ветка, а не дерево с собственными
корнями, поэтому терять связь с научным миром никак нельзя. Каковы
бы ни были олимпиадные деятели, они светятся хотя бы отраженным
светом. Нынешняя эволюция олимпиад, в том числе и в России,
предвещает мало хорошего. Если возобладает чисто спортивный
подход, то математическая олимпиада, очевидно, не сможет
конкурировать с иными соревнованиями~-- ни по зрелищности, ни по
популярности. Возникнут новые деятели, специалисты по проведению
игры типа ``завоюй красный флажок'' и они
 вытеснят старых.

Если олимпиадный деятель может отстранить ученого от участия в
подготовке олимпиады, то, наверное, и чиновнику можно заменить
олимпиадного деятеля?

Последняя олимпиадная реформа, связанная с отменой зонального
этапа, невозможностью принимать участие в городской олимпиаде, не
став победителем районной, и т.д., оказалась возможной в том числе
из-за  отсутствия в олимпиадном движении крупных ученых.
Отчуждение между научным и олимпиадным сообществами может
приводить к нежеланию последнего, чтобы в обсуждении вопросов
образования участвовали ученые и, как следствие, к сдаче  позиций
перед чиновниками. В~письме одного регионального олимпиадного
деятеля к президенту РФ говорится о привлечении энтузиастов,
работающих на местах, к обсуждению организационно-педагогических
вопросов, но не говорится о привлечении ученых.

Когда на это обстоятельство было указано автору обращения, то он
ответил, что ученые всё равно в олимпиадной деятельности не
участвуют. Такое чувство, что у сообщества пропадает инстинкт
самосохранения.  Однако, благодаря связям московских олимпиадных
деятелей с научным сообществом удалось существенно уменьшить вред
от прошедшей олимпиадной реформы.

В последнее время к функциям олимпиады добавился  способ
поступления, альтернативный ЕГЭ, что создает дополнительные
проблемы, в частности, с содержанием задач и секретностью
подготовки вариантов, c сужением круга лиц, готовящих варианты.

Имеются опасности, связанные с влиянием подготовки к международной
олимпиады на стандарты национальных олимпиад, что приводит к еще
большему сокращению тематики. Например, на международной олимпиаде
нет стереометрии. Тогда зачем ее учить? На международной олимпиаде
вариант составляется голосованием team-leader'-ов. При таком
голосовании они преследуют отнюдь не только интересы научного
содержания. Это самым плачевным образом отражается на качестве
варианта. (Хотя безусловно на международной олимпиаде бывают
красивые задачи.)

\subsection{Подход к олимпиадам, имеющий источником науку}

 {Девиз другой олимпиадной идеологии}: {\it преподавание и
олимпиады  должны отражать науку\/}. Этой идеологии следовали
всесоюзные олимпиады и старые олимпиады во многих странах (с
единственной оговоркой~-- изначально олимпиады были близки по духу
к вступительным  экзаменам).

Впоследствии выяснилось, что олимпиада~-- это, в том числе,
полигон для отработки новых тем и сюжетов (см. задачник
``Кванта'', особенно в 70-е -- 80-е годы). В современных
олимпиадах эта идеология присутствует не в чистом виде (в наиболее
чистом виде~-- в Турнире городов, особенно на его летних
конференциях \cite{lktg}, и в Московской олимпиаде, при всех их
недостатках).

Этот девиз о связи с математикой имеет конкретное преломление.

Прежде всего, химерам в олимпиадах отказано в праве на
существование. Ведь постановка задачи столь же важна, как и умение
ее решать. Сила математика, как уже говорилось, во многом зависит
от его вкуса. Совершенно необходимо иметь чутье на естественность.
Без этого человек находится вне науки. Неестественная трудная
задача на олимпиаде портит вкус и наносит огромный вред
участникам.

Хорошая олимпиадная задача получается путем оформления идей и
сюжетов из науки. Реже возникает новый сюжет в элементарной
математике.

К спорту~-- отношение утилитарное, как к средству заставить
подростка выложиться, достичь глубины, изучить технику. Приоритет
соображений, уважающих содержание, над спортивными, когда это
возможно.

Математик, занимающийся большой проблемой, ищет связанные с ней
задачи, где предполагаемые идеи решения работают в более простой
ситуации. Невозможно также придумать несколько идей сразу~-- нужны
промежуточные этапы, поэтому ценятся {\it продвижения}. Под
наличием решения, как правило, понимается наличие ``каркаса'',
~--е основных идей. Такой подход в миниатюре можно переносить на
олимпиадное творчество. Академические ценности по своей природе
неформальны, но имеют огромное значение для воспитания будущих
ученых.

К научному стилю, так или иначе, тяготеют практически все
профессиональные математики, занимавшиеся олимпиадами.

\subsection{Два подхода и человеческий фактор}

В человеческих сообществах наличие разных ценностных ориентаций
преломляется через субъективный фактор, амбиции и приводит к
``политическим'' последствиям. Понятно, что человек стремится
создать среду обитания для себя, стремится установить такую
структуру ценностей, чтобы его рейтинг поднялся. При этом от шкалы
ценностей сильно зависит выбор авторитетов и ценность тех или иных
``заслуг''. С этим связаны дополнительные проблемы при организации
олимпиадного сообщества.


В наше время происходят неприятные явления, друг друга
усиливающие. Во-первых, повторим, ослабляется влияние научного
сообщества.

Так, один филдсовский лауреат  был вовлечен в комиссии, связанные
со школьной программой. Он имел неосторожность сказать в частном
разговоре о педагогических деятелях, что доверять им образование
столь же осмысленно, как красным кхмерам вопросы демократии. Этой
устной фразы (сказанной не про конкретного человека, а вообще)
оказалось достаточно для его отстранения. Раньше в мире не было
такого влияния  невежественных деятелей в вопросах образования.

Во-вторых, ослабляется внимание ведущих ученых к вопросам
образования вообще и олимпиад в частности. Проблемы олимпиадного
движения часто аргументируют именно этим явлением. Однако эта
аргументация бывает лицемерна, поскольку главный фактор~--
вытеснение ученых.  Традиции математического образования живы, и
можно указать довольно много ученых докторского уровня, которых
можно привлечь к организации и проведению олимпиад~-- не как
``генералов'', а для практической работы.

Для создавшегося положения характерно именно сочетание этих двух
обстоятельств. Если бы ведущие ученые принимали деятельное участие
в олимпиадах, как это было 20--30 лет назад, создавшееся положение
не возникло~бы. С другой стороны, если бы не было процесса
вытеснения (об этом~-- ниже), то ряд людей докторского уровня
участвовал бы в турнирах и олимпиадах, и уровень мероприятий был
бы совсем иным.

Разумеется, без определенной преподавательской и некоторой
олимпиадной квалификации одной научной квалификации для
составителя варианта недостаточно. Действительно, существует
определенный набор вопросов, существенных при подготовке олимпиады
(наличие утешительной задачи, балансировка варианта по сложности,
наличие разных тем в варианте олимпиады). Но грамотному
математику, интересующемуся олимпиадами (а на олимпиадах
воспиталось много математиков), всему этому научиться несложно.

В учебном процессе по некоторым предметам ставится ``зачет'',
выражающий только наличие определенных умений, при этом неважно,
насколько эти умения доведены до совершенства, а по некоторым
предметам~-- ``экзамен с оценкой''. За соблюдение формата варианта
жюри можно присудить только ``зачет'' или ``незачет''. История же
ставит дифференцированную оценку только за научное содержание. У
создателей олимпиадного движения не было стремления ни к тонкой
балансировке варианта, ни к близости распределения решивших задачи
с желаемому (только чтобы оно было в разумных рамках). Эти
стремления появились и стали приводить к засилью искусственных
задач, которые (см. выше) проще придумывать и удобнее оценивать.
Такого рода ``профессионализация'' есть паразитическая часть
олимпиадной культуры. Она позволяет рационализировать отторжение
хороших задач и создает авторитет некоторым олимпиадным деятелям,
ограничивает круг лиц, занимающихся олимпиадами.

Когда за олимпиаду отвечал математик приемлемого уровня, то даже
при отсутствии олимпиадного опыта научное содержание олимпиад
выправлялось. Что касается соблюдения организационных принципов и
форматов, оно быстро приходило в норму. Так было, в частности, с
Всесоюзными и Московскими городскими олимпиадами в середине 80-х
годов. Роль сильного математика может быть нетривиальной. Бывает
очень непросто понять, как наши сильные качества или слабости
отражаются на нашем окружении. Так, только проанализировав свои
разговоры с учениками, я обнаружил, что почти все время я говорил
о технике атаки на тот или иной открытый вопрос, но не о красоте
той или иной математической теории, и осознал необходимость
дополнительных душевных усилий в этом направлении. О силе
руководителя прежде всего говорит его окружение.

\paragraph{Об олимпиадном сообществе.}
Плохо, когда педагогическое сообщество отделено от научного и
предоставлено самому себе. Часть возникающих нежелательных явлений
связана с внутренними особенностями самого элитного (по
содержанию, а не престижности) образования, которое с
неизбежностью базируется на культе успеха. Более способные
учащиеся получают больше внимания и ресурсов. (Противникам
элитного образования следует указать на их лицемерие~-- ведь в
жизни более успешный человек больше получает. В конце концов,
``ибо кто имеет, тому дано будет и приумножится, а кто не имеет, у
того  отнимется и то, что  имеет'' (Мф.~13:12).) Применительно к
математическому образованию, больше внимания получает тот, кто
лучше решает задачи.

Многие замечали, что великие тренеры часто возникают из
неудавшихся олимпиадников. Психологические комплексы, о которых
так любят говорить, несут в себе не только плохое, но и хорошее.
Зачастую, излечив человека от комплексов, мы его излечим и от
талантов и будущих достижений. Не понимая этого, гуманистические
психологи\footnote{Это их самоназвание.} делают подобные суждения
некомпетентными.

Тем не менее, могут возникать психологические проблемы, которые
следует учитывать и, по возможности, смягчать. В сообществе, в
котором доминируют неудавшиеся ученые, возникает весьма
специфическая атмосфера. Если человек уделяет большее внимание
учащемуся, более успешному в решении трудных задач, в то время как
у него самого есть проблемы с творческой реализацией (особенно
если он был вундеркиндом), он тем самым усиливает свои комплексы,
порой это становится опасным. Нужно также учитывать что психика
преподавателя зачастую приближается к подростковой.

Разумеется, отнюдь не все преподаватели и педагогические деятели
таковы. Есть вполне успешные люди, пришедшие в эту область (как
учителя, методисты или организаторы) по зову сердца, а не потому,
что наука не получилась.

Есть механизм, выдвигающий  наиболее амбициозных деятелей. Влияние
человека сильно зависит от энергии, которую он вкладывает, и, к
сожалению, от уровня его агрессии. (В этом отношении человеческое
сообщество не отличается от иной популяции животных.) Овладение
культурными ценностями, равно как и научная работа, требует очень
много энергии. Если же всего этого нет, то при том же природном
энергетическом потенциале человек будет выглядеть более ярко, чем
тот, у кого есть другая большая работа. Бурлит мелкая вода. Есть
люди с блестящим, но бесплодным интеллектом (быстрота мыслительных
операций ~-- еще не интеллект, интеллект~-- это не ум, а ум~-- не
мудрость). Однако выдвижение человека и его роль сильно зависит от
этого внешнего блеска.\footnote{Внешний блеск очень полезен в
преподавании. Но в преподавании сильным учащимся важнее научная
квалификация.} Этим объясняется парадоксальная ситуация,
сложившаяся в жюри некоторых олимпиад и турниров: ответственными
за варианты для относительно слабых учащихся (работа, связанная с
меньшими амбициями) являются более сильные математики и
методисты.\footnote{Кроме того, следует отметить,
  что студенту, пусть даже
талантливому математику и педагогу, входящему в жюри мероприятия,
иногда трудно спорить со старшим по возрасту, в особенности с тем,
с кем он взаимодействовал, будучи школьником.}

Если олимпиада~-- единственное поле реализации человека, то он
будет энергично добиваться должностей, занимаясь даже неприятной
работой, да и локтями поработает. А человека, который может
реализоваться в науке, гораздо легче прогнать~-- у него есть
другое поле деятельности. Происходит отрицательный отбор. Подобные
механизмы действия инстинкта агрессии изложены в
книгах~\cite{Dolnik,Lorenc}.

Далеко не все педагогические деятели прошлого были сильными
математиками, но они были членами научного сообщества. Есть и
сейчас деятели, высказывающие глубокие методические идеи, авторы
журнальных статей. Но зачастую задают тон люди с узкой
специализацией~-- придумывание олимпиадной задачи или создание
варианта олимпиады, профессионального с точки зрения формата
мероприятия, но с пренебрежением к содержанию.

Научный мир (особенно в точных науках), при всех своих нынешних
недостатках (см. например, книгу А.~Гротендика ``Урожаи и
посевы''), более здоровый, чем олимпиадно-педагогический. В нем
есть объективные критерии.

Но даже в более здоровом научном сообществе необходима
определенная формальная самозащита (ученые степени, звания и
т.д.), часто с издержками.
 К западному опыту
организации науки надо относиться критично, но не игнорировать.
Следует осознать причины тех или иных  решений. Там ученые
меряются, прежде всего, публикациями (с \mbox{учетом} их
количества и индекса цитируемости соответствующих журналов), а уже
затем~-- по преподаванию. Этот способ имеет много недостатков.
Зачастую очень хорошие специалисты имеют более низкий индекс
цитирования, чем некоторые ученые среднего уровня. Тем не менее,
нижний и средний уровень эти формальные параметры измеряют хорошо.
Репутация организации зависит от публикаций
сотрудников.\footnote{За рейтингами
  стоит нетривиальное управление в английском стиле.
Сравнение успехов студентов на олимпиадах (и в научной работе)
между российскими и западными университетами невыгодно Западу.
Поэтому соответствующие параметры не учитываются в подсчете
рейтингов. О системе PISA подробно рассказано в статье Д.~Малати в
журнале ``Математическое образование''. 
Понижение рейтинга отечественных журналов вызвало снижение потока
хороших статей и, как следствие, цепную реакцию понижения.}

\section{Заключение}

В создавшихся условиях особенно важны традиции, заложенные
создателями олимпиадного движения.  В частности, они  способны
затруднить, а иногда и остановить вредное реформаторство. То
обстоятельство, что у истоков олимпиад стояли великие люди, имеет
важное значение и для привлечения общественного внимания и
ресурсов (включая финансовые). Но, в то же время, это и
ответственность~-- по меньшей мере, моральная. Например, если мы
объявляем, что наша олимпиада основана таким-то ученым, или
называем мероприятие его именем, то общественность это
воспринимает как заявление в плане верности традиции, и мы тем
самым ставим вопрос: как бы относился имярек к его проведению и
стилю?\footnote{В некоторых случаях вопрос надо согласовать с
наследниками. На Западе к этому чувствительны.}

Автору близка мысль члена жюри Всероссийской олимпиады по
математике А.~С.~Голованова о роли гласности, высказанная им в
интернет-дискуссиях. О гласности в отношении уровня жюри и
достижений его членов должны заботиться организаторы
турниров.\footnote{Некоторое представление как о научных, так и о
популярных статьях, можно
получить на платных  сайтах  \texttt{http://www.zentralblatt-math.org},\\
\texttt{http://www.ams.org/mathscinet}, и бесплатных
\texttt{http://scholar.google.com/},\\
  \texttt{http://kvant.mirror1.mccme.ru/},
  \texttt{http://www.turgor.ru/lktg/index.php},\\
\texttt{http://www.mathnet.ru}).}

Избегая несправедливых обид, мы не считаем нужным называть
приверженцев того или иного подхода. Мы оцениваем не отдельных
личностей, равно как и конкретные мероприятия, но подходы и
тенденции.  Борьба разных стилей происходит внутри человека. Автор
опирается и на свой внутренний опыт. Он не может осуждать деятеля,
прервавшего свою заграничную поездку, чтобы подготовить
международника, или человека, для которого олимпиадный диплом
стоит больше публикации в научном журнале. Он сам испытывал
б\'ольшую радость, когда его задача попадала в вариант олимпиады,
чем от публикации, скажем, в Journal of Algebra или в ИАН, а
олимпиадная премия могла обрадовать больше, чем защита
диссертации. К нашим чувствам надо относиться критично. Иногда
попадаются математики, жестко следующие спортивному подходу, а с
другой стороны~-- есть педагоги, преданные
науке.\footnote{Поведение действующего математика, получившего
чрезмерно спортивное воспитание, зачастую противоречиво. Взгляды
людей, переставших заниматься наукой и оторванных от научного
сообщества, эволюционируют в сторону первого стиля. Иногда
поведение математика в отношении задач своей узкой специальности
отвечает второму подходу, а остальных задач~-- первому.} Кроме
того, сами формальные показатели нуждаются в изучении и уточнении.
Да и ситуация в каждом из конкурсов должна быть изучена подробно.

Олимпиадные функционеры тоже нужны~-- прежде всего, в технических
вопросах, а также вопросах, связанных с формированием простой
части вариантов.

Разнообразие интересов членов методкомисии улучшает атмосферу и
результат работы. Очень плохо, когда свет клином сошелся на
мероприятиях. Очень важна диверсификация интересов деятеля. Когда
интересы слишком узкие, есть очень большой риск что человек для
достижения целей способен идти на все, позабыв о цене цели и о
морали.\footnote{Такое стремление способно передаваться ученикам.
В отборочном мероприятиии на международное соревнование один
ученик дал списать другому за значительную мзду, которую тот
``отстегнул''. Списавший прошел на IMO и получил высокий приз. }

Говоря о рейтингах, кроме научных статей есть еще статьи в
достаточно престижные популярные журналы, например, ``Квант''.

Чтобы возникли новые авторы, важен малый жанр. Для написания
большой статьи надо потратить много ресурсов, да ее еще могут не
принять. Пусть и школьники, и учителя из провинции, общаясь с
редакторами научных журналов, получат навык и вкус к научным
миниатюрам.\footnote{Автор еще в 90-е годы критиковал (часто с
ущербом для себя) некоторых журнальных деятелей, ревниво
относящихся к своему делу и подавлявших других, особенно
провинциальных авторов. Многие люди, испытавшие на себе это
подавление, в последующем вели себя не лучше~-- если не хуже. Дело
не в конкретных персоналиях, а в проблемах самого педагогического
сообщества.}

Нельзя жить прошлыми заслугами. Уровень популярных журналов в
России пока выше, чем за рубежом, но это не может долго
продолжаться по инерции. Об авторитете отечественных олимпиад
также нужно заботиться, а он зависит, прежде всего, от их
содержания.

{И самое главное. Поднять научный уровень олимпиад и одновременно
спортивные достижения МОЖНО. Это проще, чем кажется. Не так сложно
научиться технике составления вариантов, да и среди действующих
математиков есть достаточное число людей, прошедших олимпиадную
школу. Организация постоянно действующего семинара по
\mbox{олимпиадным} задачам в конце 80-х -- начале 90-х годов на
мехмате привела к появлению новых олимпиадных деятелей и
достаточного запаса новых олимпиадных задач. Важно привлекать
действующих математиков, причем разнообразных специальностей,
чтобы были покрыты основные разделы науки. Тех, кого можно
привлечь к олимпиадному движению, довольно много, как в столицах,
так и на периферии.}

Обучение живой математике даже с точки зрения спортивных
достижений гораздо  эффективнее, чем натаскивание на мертвые
задачи. Этому есть примеры. Сила приводит к успехам, но
проявляется она в обличии красоты.

Данная статья призвана начать обсуждение проблем олимпиадного
движения, которое, как утверждали практически все математики и
методисты, прочитавшие эту статью, давно назрело. Необходимо
начать работу по осмыслению накопленного опыта, и это осмысление
принесет плоды. Призываем всех к дальнейшему обсуждению.

\vspace*{-\baselineskip}
\enlargethispage{2\baselineskip}

\end{document}